\documentclass[12pt]{article}
\usepackage[all]{xy}
\usepackage{amsfonts,amsmath,verbatim,oldgerm,amssymb,amscd,makeidx,latexsym}
\usepackage{theorem,calc,enumerate}
\usepackage[OT2,T1]{fontenc}
\usepackage{footnote}
\UseComputerModernTips

\newcommand{\lra}{\longrightarrow}
\newcommand{\lla}{\longleftarrow}
\newcommand{\tn}{\textnormal}

\newcommand{\os}{\overset}
\newcommand{\us}{\underset}

\newcommand{\ra}{\rightarrow}
\newcommand{\Ra}{\Rightarrow}
\newcommand{\Lra}{\Longrightarrow}

\newcommand{\sbt}{\subset}

\newcommand{\mf}{\mathfrak}

\theoremstyle{plain}
\newtheorem{theorem}{Theorem}
\newtheorem{prop}[theorem]{Proposition}
\newtheorem{lemma}[theorem]{Lemma}

\theorembodyfont{\upshape}

\newcommand{\Z}{{\mathbb Z}}
\newcommand{\Q}{\mbox{$\mathbb Q$}}
\newcommand{\R}{\mbox{$\mathbb R$}} 
\newcommand{\C}{\mbox{$\mathbb C$}} 

\newcommand{\Hom}{\operatorname{Hom}}

\newcommand{\gal}{\operatorname{Gal}}
\newcommand{\sel}{\operatorname{Sel}}

\makeindex

\def\Com{\UseComputerModernTips}

\title{Comparison of the $\mu$-invariants of an abelian variety
and its dual abelian variety}
\author{Amala Bhave\\{\small School of Physical Sciences, Jawaharlal Nehru University, New Delhi 110 067, India.}}
\date{}
\begin{document}
\maketitle

\begin{center}
Abstract
\end{center}
In this note, we compare the dual Selmer groups of an abelian variety with that
of its dual over certain large Galois field.
We give a formula which relates the generalized Iwasawa $\mu$-invariants associated
with their dual
Selmer groups under the natural isogeny.
\section{Introduction}
Let $F$ be a number field and $\bar{F}$ a fixed algebraic closure of $F$.
Let $A$ be an abelian variety defined over $F$ of dimension $d$ with no complex multiplication
over $\bar{F}$. Let $p$ be an odd prime number. By $A_{p^\infty}$, we denote the set of all $p$-th power torsion
points of $A$ and by $F_\infty$ we denote the field of definition
$F(A_{p^\infty})$ of $A_{p^\infty}$. Let $G_\infty=\gal(F_\infty/F)$. By a
well known result of Serre, if $A$ does not admit complex multiplication over
$\bar{F}$, then $G_\infty$ embeds into $GL_{2d}(\Z_p)$,

whose  
the cohomological
dimension
is the same as its manifold dimension as a $p$-adic
Lie group and hence is equal to $2d$. By $F_n$ we denote the field extension of
$F$ obtained by adjoining $p^{n+1}$-torsion points of $A$ to $F$. Let
$F^{cyc}$ be the cyclotomic $\Z_p$-extension of $F$. From the Weil pairing
it is clear that $F^{cyc}$ is contained in $F_\infty$.
Let $\Gamma=\gal(F^{cyc}/F)$.
\footnote{2010 Subject Classification. Primary $11$R$23$, $11$G$05$; Secondary $11$G$40$.}

Let $S$ be a finite set of primes of $F$ containing the all primes dividing $p$
and primes at which $A$ has bad reduction. Let
$F_S$ denotes the maximal extension of $F$ which is unramified at all the places not
contained in $S$ and at Archimedean primes.
Let $G_S=\gal(F_S/F)$. Then clearly $F_\infty$ is contained
in $F_S$.  So we have the following tower of extensions of $F$.
\Com
\begin{equation}
\begin{split}
\xymatrix{
F_S\ar@{-{}}[d]\ar@/^2.1pc/@{-{}}[ddd]^{G_S}\\
F_\infty\ar@{-{}}[d]^H\ar@/_1.8pc/@{-{}}[dd]_{G_\infty}\\
F^{cyc}\ar@{-{}}[d]^{\Gamma}\\
F
}
\end{split}
\end{equation}

For any extension $L$ of $F$ contained in 
$F_S$, let $G_S(L)$ denote the Galois group
$\gal(F_S/L)$. Then for any extension $L$ of $F$
contained in $F_\infty$ the
$p^\infty$-Selmer group of $A$ over $L$ denoted by $\sel_{p^\infty}(A/L)$ is
defined by the exact sequence
\begin{equation}\label{selmer}
0\lra\sel_{p^\infty}(A/L)\lra \tn{H}^1(G_S(L),A_{p^\infty})\os{\lambda_L}{\lra}
\us{v\in S}{\oplus}J_v(A/L)
\end{equation}
where $J_v(A/L)=\us{\lra}{\lim}\us{w|v}{\oplus}\tn{H}^1(K_w,A)(p)$, $K$ runs over
the finite extension of $F$ contained in $L$ and the direct limit is taken with
respect to the restriction maps.

For any profinite group $G$, 
let $\Lambda(G):=\us{\lla}{\lim}~\Z_p[G/U]$, ($U$ being open normal
subgroup of $G$) be the Iwasawa algebra of $G$ which is a complete ring.
The Galois group $G_\infty$ acts on $\sel_{p^\infty}(A/F_\infty)$ continuously
and the action is extended continuously so that it becomes a discrete module
over the Iwasawa algebra $\Lambda(G_\infty)$. We consider
its Pontrjagin dual
${\mf X}(A/F_\infty)=\widehat{\sel_{p^\infty}(A/F_\infty)}$
which is a compact
$\Lambda(G_\infty)$-module. Similarly $\mf{X}(A/F^{cyc})$ is a
$\Lambda(\Gamma)$-module. It is known (\cite{1}) that $\mf{X}(A/F^{cyc})$ is a finitely
generated $\Lambda(\Gamma)$-module and that ${\mf X}(A/F_\infty)$ is a
finitely generated $\Lambda(G_\infty)$-module.

The following is conjectured.\\
{\bf Conjecture 1:} Let $F^{cyc}\sbt L\sbt F_S$. Then \\
i) the map $\lambda_L$ in the exact sequence \eqref{selmer} is surjective. \\
ii) $\tn{H}^2(G_S(L),A_{p^\infty})=0$.

It is already known (Theorem 2.10 in \cite{2}) that
$\tn{H}^2(G_S(F_\infty),A_{p^\infty})=0$ but this is still
only a conjecture for an arbitrary $L$. \\
{\bf Conjecture 2:} $\mf{X}(A/F_\infty)$ is a torsion
$\Lambda(G_\infty)$-module.

We assume both these conjectures. In fact conjecture 1 is equivalent to
conjecture 2 if $A$ has potentially ordinary reduction at all the
places dividing $p$.
If $G_\infty$ is pro-$p$ and if $M$ is a finitely generated
$\Lambda(G_\infty)$-module, then (\cite{4}) we have

\begin{equation}
\mu(M)=\mu_{G_\infty}(M)=\log_p(\chi(G_\infty,\widehat{M(p)})).
\end{equation} 


With the above notation, we state the formula for the $\mu$-invariants
of isogenous abelian varieties which has been proved by Susan Howson (\cite{4}) in
the case of elliptic curves.
\begin{prop} \label{abelianvarformula}
Let $A_1$ and $A_2$ be non-CM abelian varieties defined over $F$ and let
$f:A_1\ra A_2$ be an isogeny defined over $F$. We assume that $F_\infty=
F(A_{1,p^\infty})=F(A_{2,p^\infty})$ be a pro-$p$-extension of $F$.
Let $C=(ker~f)(p)$ and $\tilde{C_v}$ be the image of $C$ under the
reduction map at $v: A_1\ra \tilde{A_1}$. Then as
$\Lambda(G_\infty)$-modules
\begin{eqnarray}\label{formula1}
\mu(\mf{X}(A_1/F_\infty))-\mu(\mf{X}(A_2/F_\infty))~~~~~~~~~~~~~~~~~~~~~
~~~~~~~~~~~~~~~~~~~~~~~~~~~\nonumber\\
=\sum_{v|\infty}\log_p(\# C(F_v))-[F:\Q]\log_p(\# C)+
\sum_{v|p}\log_p(|\#\tilde{C_v}|_v)
\end{eqnarray}
where $|\cdot|_v$ is the $v$-adic norm.
\end{prop}
\noindent Proof.~
The argument is exactly the same as that in the $GL_2$ case (\cite{4}) in the case of an isogeny
between two abelian varieties with $F_\infty=F(A_{i,p^\infty})$. Hence we omit it.\\

\hfill$\square$
\\

We remark that Proposition \ref{abelianvarformula} has been proved for the
cyclotomic extension by Perrin-Riou (\cite{7} and
P. Schneider(\cite{PSchn}). 
\section{Comparing the abelian variety with its dual}
In this section, we prove that the $\mu$-invariants of an abelian variety
and that of its dual abelian variety are the same. 

Suppose $A_1=A$ and $A_2=A^t$, the dual abelian variety.
We always have (\cite{6}) an isogeny $f:A_1\lra A_2$.

The main ingredients of the proof are the following lemmas. For a finite
module $M$, we define $M^D:=\Hom(M,\mu_{p^\infty})$. Recall that
$C$ is the kernel of the isogeny $f:A\ra A^t$. For a prime
$v$, the kernel of the reduction of $C$ modulo $v$ is denoted by
${\cal C}_v$.
\begin{lemma} \label{primesabovep}
Let $v$ be a prime of $F$ dividing $p$. Maintaining the above notation
and hypotheses on the variety $A$, we have
$$\#{\cal C}_v=\#\tilde{C_v}=\#\tilde{C^t_v})=\#{\cal C}^t.$$
\end{lemma}
\noindent Proof.~
Let ${\cal A}_v$ be the kernel of the reduction of $A$ modulo $v$. Then we
have the following commutative diagram:
\Com 
$$
\xymatrix{
&0 \ar[d] &0\ar[d] &0\ar[d]\\
0 \ar[r] & {\cal C} \ar[r]\ar[d]  &{\cal A}_{v,p^\infty} \ar[r]\ar[d]
 & {\cal A}^t_{v,\infty} \ar[r]\ar[d] & 0\\
0 \ar[r] & C \ar[r]\ar[d]  &A_{p^\infty} \ar[r]\ar[d]
 & A^t_{p^\infty} \ar[r]\ar[d] & 0\\
0 \ar[r] & \tilde{C_v} \ar[r]\ar[d]  &\tilde{A}_{v,p^\infty} \ar[r]\ar[d]
 & \tilde{A^t}_{v,\infty} \ar[r]\ar[d] & 0\\
&0&0&0\label{diagram2}
}
$$ By Milne \cite{6} we know that $C^t=\Hom(C,\mu_{p^\infty})$ which gives
perfect Weil pairing $C \times C^t \lra \mu_{p^\infty}$. Further the exact
annihilator of ${\cal A}_{p^n}$ is ${\cal A}^t_{p^n}$ and therefore we see that
${\cal A}_{p^n}$ is Pontrjagin dual to $\tilde{A^t}_{p^n}(-1)$. Using this in
the  pairing $C \times C^t \lra \mu_{p^\infty}$, we see that the Pontrjagin
dual of ${\cal C}_v$ is $\tilde{C^t}_v(-1)$. This implies that
$\#{\cal C}=\tilde{C^t}_v$. Similarly, $\#{\cal C}^t=\tilde{C_v}$. We now show
that $\#{\cal C}_v=\#\tilde{C_v}$.

The snake lemma induced by multiplication by $p^n$ for the top row in 
the above diagram
gives an exact sequence
$$0\lra{\cal C}_{v,p^n}\lra{\cal A}_{v,p^n}\lra{\cal A}^t_{v,p^n}\lra
\frac{{\cal C}_v}{p^n{\cal C}_v}\lra0$$ since ${\cal A}_{v,p^\infty}$ is
divisible. Taking inverse limit over $n$, we get an exact sequence
$$0\lra T_p{\cal A}_v\lra T_p{\cal A}^t_v\lra {\cal C}_v\lra0$$
since ${\cal C}_v$ is a finite $p$-primary group and hence
$\us{\lla}{\lim}~{\cal C}_{v,p^n}=0$.

Next recall that ${\cal A}_{v,p^n}$ is Pontrjagin dual to
$\tilde{\cal A}^t_{v,p^n}(-1)$. This implies, $$0\lra
\widehat{\tilde{\cal A}^t_{v,p^n}(-1)}\lra
\widehat{\tilde{\cal A}_{v,p^n}(-1)}\lra{\cal C}_v\lra0.$$ Taking the
Pontrjagin dual along this exact sequence, we have
$$0\lra\hat{\cal C}_v\lra\tilde{\cal A}_{v,p^n}(-1)\lra
\tilde{\cal A}_{v,p^n}^t(-1)\lra0.$$ But by the bottom row of diagram
\eqref{diagram2}, $\#\hat{\cal C}_v=\#\tilde{C_v}(-1)$ and hence
$\#{\cal C}_v=\#\tilde{C_v}$.
Hence the lemma.\\

\hfill$\square$
\begin{lemma}\label{infiniteprimes}
If $v|\infty$ and $G_v=\gal(\bar{F}_v/F_v)$, then $\#C(F_v)=\chi(F_v,C)$
and $\#C^t(F_v)=\chi(F_v,C^t)$.
\end{lemma}
\noindent Proof.~
First let $v$ be a complex prime so that $F_v=\C$. Then
$G_v=\{e\}$ is the trivial group. This implies
$\tn{H}^1(G_v,C)=\tn{H}^2(G_v,C)=0$ and hence $\chi(G_v,C)=\#\tn{H}^0(G_v,C)=\#C$.
Similarly, $\chi(G_v,C^t)=\#\tn{H}^0(G_v,C^t)=\#C^t$.

Now let $v$ be a real prime so that $F_v=\R$ and $\bar{F}_v=\C$ with
$G_v=\Z/2\Z$. As $G_v$ is a finite cyclic group, its Herbrand quotient is
$1$. Therefore $\chi(G_v,C)=\#\tn{H}^0(G_v,C)=\#C(F_v)$. Similar argument for
the dual abelian variety gives $\chi(G_v,C^t)=\#C^t(F_v)$.\\

\hfill$\square$
\begin{prop}\label{samemuinvariant}
With same notation and assumptions as bove, 
$$\mu({\mf X}(A/F_\infty))=\mu({\mf X}(A^t/F_\infty)).$$
\end{prop}
\noindent Proof.~
By Proposision \ref{abelianvarformula}, we have the formula
$$\mu({\mf X}(A/F_\infty))-\mu({\mf X}(A^t/F_\infty))=\sum_{v|\infty}
\log_p(\#C(F_v))-[F:\Q]\log_p(\#C)$$
$$~~~~~~~~~~~+\sum_{v|p}\log_p|\#\tilde{C_v}|_v;$$ and similarly by
interchanging the roles of $A$ and $A^t$ we have another formula
$$\mu({\mf X}(A^t/F_\infty))-\mu({\mf X}(A/F_\infty))=\sum_{v|\infty}
\log_p(\#C^t(F_v))-[F:\Q]\log_p(\#C^t)$$
$$~~~~~~~~~~~+\sum_{v|p}\log_p|\#\tilde{C_v}^t|_v.$$
By Lemma \ref{primesabovep}, $\#\tilde{C_v}^t=\tilde{C_v}$ for primes lying
above $p$ and by Lemma \ref{infiniteprimes}, we have equality
$\prod_{v|\infty} \#C(F_v)= \prod_{v|\infty}\#C^t(F_v)$. These two together
imply that the
left hand side expressions of the above two formulae are equal. But clearly,
they are negatives of each other. Hence left hand side expression is zero.
That is $\mu({\mf X}(A/F_\infty))=\mu({\mf X}(A^t/F_\infty))$. Hence the
proposition.\\

\hfill$\square$
\\

In the cyclotomic case, Proposition \ref{samemuinvariant} has been
proved by K. Matsuno (\cite{5}). Our methods give a different proof of his
result which carries over to the $GL_2$-case. It is not a priori clear how to
generalize Matsuno's proof.

%

A\small{MALA} B\small{HAVE}\\
School of Physical Sciences\\
Jawaharlal Nehru University\\
New Delhi 110 067\\
India\\
E-mail: abhave@mail.jnu.ac.in
\end{document}